\documentclass[MSNbibl,number,citesort,dvips]{arxbj}
\usepackage{upgreek}
\usepackage{dcolumn}

\aid{0}
\volume{19}
\issue{2}
\pubyear{2013}
\firstpage{409}
\lastpage{425}
\doi{10.3150/11-BEJ411}

\makeatletter
\newcolumntype{d}[1]{D{.}{.}{#1}}%lygiavimas pagal taska
\def\R{\mathbb{R}}
\def\N{\mathbb{N}}
\def\K{\mathcal{K}}

\newremark{ex}{Example}[section]
\newproclaim{con}[ex]{Condition}
\newtheorem{theo}[ex]{Theorem}
\makeatother

\begin{document}
\begin{frontmatter}

\title{Statistical inference for discrete-time samples from
affine stochastic delay differential equations}
\runtitle{Statistical inference for stochastic delay differential equations}

\begin{aug}
%%%% inicialai - be tarpu
\author[1]{\fnms{Uwe} \snm{K\"{u}chler}\thanksref{1}\ead[label=e1]{kuechler@mathematik.hu-berlin.de}} \and
\author[2]{\fnms{Michael} \snm{S{\O}rensen}\corref{}\thanksref{2}\ead[label=e2]{michael@math.ku.dk}}
\runauthor{U. K\"{u}chler and M. S{\o}rensen} %% auto
\address[1]{Institut f\"ur Mathematik, Humboldt-Universit\"at zu Berlin,
Unter den Linden 6, D-10099 Berlin, Germany. \printead{e1}}
\address[2]{Department of Mathematical Sciences, University of Copenhagen,
Universitetsparken 5, \mbox{DK-2100} Copenhagen {\O}, Denmark. \printead{e2}}
\end{aug}

% HISTORY:
\received{\smonth{4} \syear{2009}}

% ABSTRACT
%
\begin{abstract}
Statistical inference for discrete time observations of an affine stochastic
delay differential equation is considered. The main focus is on maximum
pseudo-likelihood estimators, which are easy to calculate in practice.
A more general class of prediction-based estimating functions
is investigated as well. In particular, the optimal prediction-based
estimating function and the asymptotic properties of the
estimators are derived. The maximum pseudo-likelihood estimator is a
particular case, and an expression is found for the efficiency loss
when using the maximum pseudo-likelihood estimator, rather than the
computationally more involved optimal prediction-based estimator.
The distribution of the pseudo-likelihood estimator is investigated in a
simulation study. Two examples of affine stochastic delay equation are
considered in detail.
\end{abstract}

% KEYWORDS
%
\begin{keyword}
\kwd{asymptotic normality}
\kwd{composite likelihood}
\kwd{consistency}
\kwd{discrete time observation of continuous-time models}
\kwd{prediction-based estimating functions}
\kwd{pseudo-likelihood}
\kwd{stochastic delay differential equation}
\end{keyword}

\end{frontmatter}
%
%s1 ###
\section{Introduction}
In the last decade, statistical inference for stochastic delay
differential equations (SDDEs) has been studied from various
viewpoints. Early work on maximum likelihood estimation was done by
{K\"uchler} and Mensch~\cite{kuechmen}.
Gushchin and {K\"uchler}~\cite{gishkuech99} and {K\"uchler} and
Kutoyants~\cite{kuechkut}
determined the non-standard asymptotic properties of the maximum
likelihood estimator for SDDEs, and {K\"uchler} and Vasil'jev~\cite{kuechvas05b}
constructed sequential
procedures with a given accuracy in the $L_2$ sense. Nonparametric
estimators for affine SDDEs were investigated by {Rei\ss}~\cite{reiss1} and
{Rei\ss}~\cite{reiss2}. All these studies were concerned with continuous
observation of the solution process.

As opposed to the situation for ordinary stochastic differential
equations, observations at discrete time points have been
little studied for SDDEs. {Rei\ss}~\cite{reiss} studied nonparametric
estimation. {K\"uchler} and {S\o rensen}~\cite{kuso2} proposed a simple
estimator for
the parameters $\alpha_k$ in the particular type of SDDE given by
(\ref{affinN}) below. This estimator is biased, however, and can only
be expected to work well for high-frequency observations. In this
paper we report a first attempt at investigating parametric
inference for affine stochastic delay equations of the general type
(\ref{geneqn}) observed at
discrete time points. We propose a pseudo-likelihood function and
study it in the framework of prediction-based estimating functions.
Applying the methods proposed here in practice often requires the
ability to simulate solutions of SDDEs. This problem has
been studied by, among other, {K\"uchler} and Platen \cite
{kuechplaten00} and
Buckwar~\cite{buckwar}. A practical application of one of the simplest SDDEs,
discussed in
Example~\ref{affinex2} below, was provided by {K\"uchler} and Platen
\cite{kuechlerplaten07}.

We consider the model given by the stochastic differential equation
%
%e1 ###
%
\begin{equation}
\label{geneqn}
\mathrm{d}X(t) = \biggl( \int_{-r}^0 X(t+s)a_\alpha(\mathrm{d}s)\biggr )\,\mathrm{d}t
+ \sigma\,\mathrm{d}W(t),
\end{equation}
where $a_\alpha$ is a measure on $[-r,0]$ ($0 \leq r < \infty$)
such that (\ref{geneqn}) has a unique stationary solution (for a
suitable given initial condition). Conditions under which
(\ref{geneqn}) has a unique stationary solution were given by
Gushchin and {K\"uchler}~\cite{gishkuech}. By Theorem 3.1 in this
paper, the stationary solution
is a Gaussian process. We assume that the measure $a_\alpha$ depends
on a parameter $\alpha$. The parameter about which inference is to be
drawn is $(\alpha,\sigma)$ or $(\alpha,\sigma, r)$, ($\sigma
,r>0$). As
usual, we denote the parameter space by $\Theta\subseteq\R^p$. The
process $W$ is a Wiener process. The initial
condition is that the distribution of $\{ X(s) \mid s \in[-r,0]
\}$ is the stationary distribution, which always has expectation 0.
The data are observations at discrete time points: $X(\Delta), X(2
\Delta), \ldots, X(n \Delta)$.

An interesting particular case of (\ref{geneqn}) is
%
%e2 ###
%
\begin{equation}
\label{affinN}
\mathrm{d}X(t) = \sum_{k=1}^N \alpha_k X(t-r_k)\, \mathrm{d}t + \sigma \,\mathrm{d}W(t).
\end{equation}
Here the measure $a_\alpha$ is concentrated in the discrete points
$-r_1, \ldots, -r_N$, ($r_i \geq0$). The vector $(r_1, \ldots,
r_N)$ may be among the parameters to be estimated. The particular case
where $N=2$ and $r_1=0$ is considered in detail in Example~\ref{affinex2}.

In Section~\ref{sec2} we discuss how to calculate the likelihood function for
discrete time observations, and propose a pseudo-likelihood function
that closely approximates the likelihood function and is considerably
easier to calculate. We consider two examples in detail. In Section
\ref{sec3} we present prediction-based estimating functions for affine stochastic
delay equations, find the optimal estimating function in this class,
and show that the pseudo-likelihood estimator is a
particular case of a prediction-based estimator. The
prediction-based estimating functions provide a good framework for
discussing the asymptotics of the pseudo-likelihood estimator and in
particular the efficiency loss compared with the optimal
prediction-based estimating function. We do this in Section~\ref{sec4},
specifying conditions ensuring consistency and asymptotic normality.
Finally, in Section~\ref{simstudy} we present a simulation study
of properties of the pseudo-likelihood
estimator.
%
%s2 ###
\section{The likelihood and the pseudo-likelihood function}\label{sec2}
Because the stationary solution to (\ref{geneqn}) is a zero-mean
Gaussian process, Gushchin and {K\"uchler}~\cite{gishkuech}, the data
are in fact a Gaussian
time series with expectation 0. Therefore,
in principle the likelihood function can be calculated if we can determine,
analytically or numerically, the autocovariances
%
%e3 ###
%
\begin{equation}
\label{genK}
K_\theta(t) = E_\theta(X(0)X(t)), \qquad t \geq0.
\end{equation}
The autocovariance function, $K_\theta(t)$, satisfies the differential
equation
%
%e4 ###
%
\begin{equation}
\label{YW}
\partial_t K_\theta(t) = \int_{-r}^0 K_\theta(t+s) a_\theta(\mathrm{d}s), \qquad t
\geq0,
\end{equation}
with $\partial_t K_\theta(0+) = - \frac12 \sigma^2$, provided that
we define $K_\theta(-t) = K_\theta(t)$ for $t \geq0$ (see
Gushchin and {K\"uchler}~\cite{gishkuech03}). The condition $\partial_t
K_\theta(0+) = -
\frac12 \sigma^2$ also can be written in the form
\[
2 \int_{-r}^0 K_\theta(s) a_\theta(\mathrm{d}s) = -\sigma^2.
\]
Equation (\ref{YW}) is a continuous-time analogue of the
Yule--Walker equation known from time-series analysis, and hereinafter we
refer to (\ref{YW}) as the delay Yule--Walker equation of
(\ref{genK}). In general, this
equation must be solved numerically, but we consider two
particular examples where it can be solved explicitly.

To calculate the \textit{likelihood function},
define for every $\ell= 1, \ldots, n$ the $\ell$-dimensional
vector $\kappa_\ell(\theta) = (K_\theta(\Delta), \ldots,
K_\theta(\ell\Delta))^T$,
and the $\ell\times\ell$-matrix $\K_\ell(\theta) = \{ K_\theta
((i-j)\Delta)
\}_{i,j = 1, \ldots, \ell}$. Here and later $T$ denotes transposition of
vectors and matrices. The matrix $\K_\ell(\theta)$ is the covariance
matrix of the vector of the first $\ell$ observations $X(\Delta),
\ldots, X(\ell\Delta)$.

The conditional distribution of the observation
$X((i+1) \Delta)$ given the previous
observations $X(\Delta), \ldots, X( i \Delta)$ is the Gaussian
distribution with expectation $\phi_i(\theta)^T X_{i:1}$ and variance
$v_i(\theta)$, where $\phi_i(\theta)$ is the $i$-dimensional vector
given by $\phi_i(\theta) = \K_i(\theta)^{-1} \kappa_i(\theta) $,
$v_i(\theta) = K_\theta(0) - \kappa_i(\theta)^T \K_i(\theta)^{-1}
\kappa_i(\theta)$, and $X_{i:j} = (X(i \Delta), \ldots, X(j \Delta
))^T$, $i > j \geq1$. The vector $\phi_i(\theta) =
(\phi_{i,1}(\theta) , \ldots, \phi_{i,i}(\theta))^T$
and the conditional variance $v_i(\theta)$ can be found using the
Durbin--Levinson algorithm (see, e.g., page 169 in Brockwell and Davis
\cite{BrockwellDavis}).
Specifically,
$\phi_{1,1}(\theta) = K_\theta(\Delta)/K_\theta(0)$ and
$v_0(\theta) =
K_\theta(0)$, whereas
\begin{eqnarray*}
\phi_{i,i}(\theta) &=& \Biggl( K_\theta(i \Delta) - \sum_{j=1}^{i-1}
\phi_{(i-1),j}(\theta) K_\theta\bigl((i-j)\Delta\bigr) \Biggr) v_{i-1}(\theta)^{-1},\\
\pmatrix{\phi_{i,1}(\theta)\cr \vdots\cr \phi_{i,i-1}
(\theta)}&=&\pmatrix{\phi_{i-1,1}(\theta)\cr
\vdots\cr \phi_{i-1,i-1}(\theta)}- \phi_{i,i}(\theta)\pmatrix{\phi_{i-1,i-1}(\theta)
\cr \vdots\cr \phi_{i-1,1}(\theta)}
\end{eqnarray*}
and
\[
v_i(\theta) = v_{i-1}(\theta) \bigl( 1 - \phi_{i,i}(\theta)^2 \bigr).
\]

The likelihood function based on the data $X(\Delta), \ldots, X( n
\Delta)$
is
%
%e5 ###
%
\begin{eqnarray}
\label{lik}
L_n(\theta) &=& \frac{1}{\sqrt{2 \uppi v_0(\theta)}}
\exp\biggl( - \frac{1}{2 v_0(\theta)} X(\Delta)^2 \biggr)\nonumber
\\[-8pt]
\\[-8pt]
&&{}
\times\prod_{i=1}^{n-1} \biggl[ \frac{1}{\sqrt{2 \uppi v_i(\theta)}}
\exp\biggl( - \frac{1}{2 v_i(\theta)} \bigl( X\bigl((i+1) \Delta\bigr) -
\phi_i(\theta)^T X_{i:1} \bigr)^2 \biggr) \biggr]. \nonumber
\end{eqnarray}
Calculation of this function quickly becomes very time-consuming as the
sample size $n$ increases. In particular, $\phi_i(\theta)$ and
$v_i(\theta)$ must be calculated for every observation
time point. However, the autocovariances $K_\theta(i \Delta)$ decrease
exponentially with $i$ (see Diekmann et~al.~\cite{diekmann}, page 34).
Using the
Durbin--Levinson algorithm, it is readily apparent that this implies
that the
quantities $\phi_{i,j}(\theta)$ decrease exponentially with $j$. Thus
the conditional distribution of $X((i+1) \Delta)$ given
$X(\Delta), \ldots, X( i \Delta)$ depends only very little on observations
in the distant past.

Therefore, we propose using instead a \textit{pseudo-likelihood
function} obtained by replacing in the likelihood function the conditional
density of $X((i+1) \Delta)$ given $X(\Delta), \ldots, X( i \Delta)$
with the conditional density of $X((i+1) \Delta)$ given
$X((i+1-k)\Delta),
\ldots, X( i \Delta)$, where $k$ typically is relatively small.
This pseudo-likelihood function was proposed by H. {S\o rensen}~\cite{helle2001}
in connection with stochastic volatility models, but the idea is widely
applicable. The pseudo-likelihood is given by
%
%e6 ###
%
\begin{equation}
\label{pslik}
\tilde{L}_n(\theta) =
\prod_{i=k}^{n-1} \biggl[ \frac{1}{\sqrt{2 \uppi v_k(\theta)}}
\exp\biggl( - \frac{1}{2 v_k(\theta)} \bigl( X\bigl((i+1) \Delta\bigr) -
\phi_k(\theta)^T X_{i:i+1-k} \bigr)^2 \biggr) \biggr].
\end{equation}
We have not included the density of $X_{k:1}$. Note that the
computational gain is large because we calculate (\ref{pslik}) using
the same values of $\phi_k(\theta)$ and $v_k(\theta)$ for all
observation time points. Thus, these quantities must be calculated
only once for every value of $\tilde{L}_n(\theta)$. We call the
number $k$ the \textit{depth} of the pseudo-likelihood function. We
consider the
influence of $k$ on the quality of the estimators in the
simulation study reported in Section~\ref{simstudy}. As would be
expected, the
quality increases with increasing depth. For the model
considered in Section~\ref{simstudy}, the present study indicates that the
bias and the variance of the estimators do not depend much on the
depth when $k$ is larger than 3--5 times $r$.
\begin{ex}
\label{affinex2}
Consider the equation
%
%e7 ###
%
\begin{equation}
\label{affin}
\mathrm{d}X(t) = [a X(t) + b X(t-r) ]\,\mathrm{d}t + \sigma \,\mathrm{d}W(t),
\end{equation}
where $r>0$, $\sigma> 0$. This is a particular case of the model
(\ref{affinN}). The real parameters $a$ and $b$ are
chosen such that a stationary solution of (\ref{affin}) exists.
This is the case exactly when $a < r^{-1}$ and $-a/\cos(\xi(ar)) < b
< -a$ if $a \neq0$, and when $-\uppi/2 < br < 0$ if $a=0$. Here the function
$\xi(u) \in(0, \uppi)$ is the root of $\xi(u) = u \tan(\xi(u))$ if
$u \neq0$, and $\xi(0) = \uppi/2$. The stationary solution is unique if
it exists. Details of this have been provided by {K\"uchler} and Mensch
\cite{kuechmen}, who
explicitly found the covariance function of the stationary solution
by solving the Yule--Walker delay differential
equation (\ref{YW}),
\[
\partial_t K_\theta(t) = a K_\theta(t) + b K_\theta(t-r), \qquad t \geq0.
\]
They found that
%
%e8 ###
%
\begin{equation}\label{affinK0}
K_\theta(0)  =\cases{\displaystyle\frac{ \sigma^2 ( b \sinh(\lambda(a,b)r ) -
\lambda(a,b))}{ 2 \lambda(a,b)[a + b \cosh(\lambda(a,b)r)] } & \quad $\mbox{when } |b| < - a$,\vspace*{+2pt}\cr
\displaystyle\sigma^2 (b r -1)/(4 b) & \quad $\mbox{when } b = a$,\vspace*{+2pt}\cr
\displaystyle\frac{ \sigma^2 (b \sin(\lambda(a,b)r ) -
\lambda(a,b))}{ 2 \lambda(a,b)[a + b \cos(\lambda(a,b)r)] } &\quad$\mbox{when } b < - |a|$,}
\end{equation}
where $\lambda(a,b) = \sqrt{|a^2 - b^2|}$, and that for
$t \in[0,r]$ the covariance function is
%
%e9 ###
%
%
\begin{equation}
\label{affinKt0}
\hspace*{-10pt}
K_\theta(t)  =\cases{K_\theta(0) \cosh( \lambda(a,b) t) - \sigma^2 (2 \lambda
(a,b))^{-1} \sinh
( \lambda(a,b) t) & \quad $\mbox{when } |b| < - a$,\vspace*{+2pt}\cr
K_\theta(0) - \frac12 t \sigma^2 & \quad $\mbox{when } b = a $,\vspace*{+2pt}\cr
K_\theta(0) \cos( \lambda(a,b) t) - \sigma^2 (2 \lambda(a,b))^{-1}
\sin
( \lambda(a,b) t)  &\quad$\mbox{when } b < - |a|$.}
\end{equation}
Because $K_\theta(t)$ is known in $[0,r]$, the Yule--Walker
equation becomes an ordinary differential equation for
$K_\theta(t)$ in $[r,2r]$, which can be easily solved. Similarly, for
$t > r$, the autocovariance function $K_\theta(t)$ is
given by
%
%e10 ###
%
\begin{equation}
\label{K2r}
 K_\theta(t) = b \int_{nr}^t \mathrm{e}^{a(t-s)}K_\theta(s-r)\,\mathrm{d}s +
\mathrm{e}^{a(t-nr)}K_\theta(nr),\qquad t \in[nr,(n+1)r], n \in\N.
\end{equation}
Thus $K_\theta(t)$ can be determined iteratively in each of the
intervals $t \in[nr,(n+1)r]$, $n \in\N$.
Note that the covariance function depends on $\sigma$ and $r$ in a
simple and smooth way, so that these parameters also can be estimated by
maximizing the pseudo-likelihood function (\ref{pslik}).

For $b=0$, the model (\ref{affin}) is the Ornstein--Uhlenbeck process, for
which (\ref{affinK0}) and (\ref{affinKt0}) simplifies to the
well-known result $K_\theta(t) = - (\sigma^2/(2 a)) \mathrm{e}^{at}$ $(t
\geq0)$ in
the stationary case $a < 0$. For $a=0$, we obtain the model
%
%e11 ###
%
\begin{equation}
\label{affin1}
\mathrm{d}X(t) = b X(t-r)\,\mathrm{d}t + \sigma \,\mathrm{d}W(t).
\end{equation}
This process is stationary if and only if $br \in(-\uppi/2, 0)$, and in
this case, by
(\ref{affinK0}) and (\ref{affinKt0}), the autocovariance function
is given by
%
%e12 ###
%
\begin{equation}
\label{affin1Kt}
K_\theta(t) = - \frac{\sigma^2}{2b}\biggl ( \frac{1-\sin(br)}{\cos(br)}
\cos(bt) + \sin(bt)\biggr )
\end{equation}
when $t \in[0,r]$. By (\ref{K2r}), we find that
%
%e13 ###
%
\begin{eqnarray}
\label{affin1Kt2r}
K_\theta(t) =
- \frac{\sigma^2}{2b} \bigl[ 2 + \cos(bt) \bigl\{ \bigl(\tan(bt) - \tan(br)\bigr)\bigl(1 -
2 \sin(br)\bigr) -1/\cos(br)\bigr \} \bigr]
\end{eqnarray}
for $t \in[r,2r]$.\vadjust{\goodbreak}
\end{ex}
\begin{ex}
\label{affinex4}
Consider the equation
%
%e14 ###
%
\begin{equation}
\label{exp}
\mathrm{d}X(t) = -b \biggl( \int_{-r}^0 X(t+s) \mathrm{e}^{as} \,\mathrm{d}s \biggr) \,\mathrm{d}t +
 \sigma\,\mathrm{d}W(t),\vspace*{-2pt}
\end{equation}
where $r>0$, $\sigma> 0$. The set of values of the parameters $a$ and $b$
for which a unique stationary solution of (\ref{exp}) exists was studied
by {Rei\ss}~\cite{reiss}. This set is rather complicated and irregular;
for instance,
it is not convex. However, it contains the region
$ \{ (a,b) \mid a \geq0, b > 0, b(1 + \mathrm{e}^{-ar}) < \max(\uppi^2/r^2, a^2
(\mathrm{e}^{ar} -1)^2) \}$. For $a = 0$, corresponding to a uniform delay measure,
a stationary solution exists exactly when $0 < b < \frac12 \uppi
^2/r^2$. When
$r = \infty$, the situation is much simpler. In that case, a stationary
solution exists for all $a>0$ and $b>0$.\looseness=-1

When $a=0$ (and $r$ is finite),
\[
K_\theta(t) = \frac{\sigma^2 \sin( r \sqrt{2 b}(1/2 - t) )}
{2 r \sqrt{2 b} \cos( r \sqrt{b/2} )} + \frac{\sigma^2}{2
br^2},\qquad 0 \leq t \leq r.
\]
For $a>0$, an explicit expression for $K_\theta(t)$ involving trigonometric
functions exists as well (see
{Rei\ss}~\cite{reiss}, page 41), but it is somewhat complicated, and thus
we omit
it here.\vspace*{-3pt}
\end{ex}
%
%s3 ###
\section{Prediction-based estimating functions}\label{sec3}
In this section, we discuss the pseudo-likelihood estimator in the framework
of prediction-based estimating functions. This class of estimating functions
was introduced by {S\o rensen}~\cite{ms2000} as a generalization of the
martingale
estimating functions that is also applicable to non-Markovian
processes such as solutions to stochastic delay differential
equations. Applications of the methodology to observations of
integrated diffusion processes and sums of diffusions have been
described by Ditlevsen and {S\o rensen}~\cite{ditlevsensor} and Forman
and {S\o rensen}~\cite{forman}. An up-to-date
review of the
theory of prediction-based estimating functions has been provided by
{S\o rensen}~\cite{ms2011}.

We show that the pseudo-likelihood estimator is a
prediction-based estimator, and find the optimal prediction-based
estimating function, which turns out to be different from the pseudo-score
function. Optimality is in the sense of Godambe and Heyde \cite
{godambeheyde} (see
Heyde~\cite{heyde97}). We impose the following condition that is satisfied
for the models considered in Examples~\ref{affinex2}
and~\ref{affinex4}.\vspace*{-2pt}
\begin{con}
\label{conA0}
The function $K_\theta(t)$ is continuously differentiable with respect to
$\theta$.
\end{con}

Under this assumption, we find the following expression
for the pseudo-score function:
%
%e15 ###
%
\begin{eqnarray}
\label{psscore}
\partial_\theta\tilde{\ell}_n(\theta) &:= &\partial_\theta\log(
\tilde{L}_n (\theta)) \nonumber\\[-2pt]
&=& \sum_{i=k}^{n-1}
\frac{\partial_\theta\phi_k(\theta)^T X_{i:i+1-k}}{v_k(\theta)}
\bigl( X\bigl((i+1) \Delta\bigr) - \phi_k(\theta)^T X_{i:i+1-k} \bigr)
\\[-2pt]
&&{}+ \frac{\partial_\theta v_k(\theta)}{2 v_k(\theta)^2} \sum_{i=k}^{n-1}
\bigl[ \bigl( X\bigl((i+1) \Delta\bigr) - \phi_k(\theta)^T X_{i:i+1-k} \bigr)^2
- v_k(\theta)\bigr ]. \nonumber\vadjust{\goodbreak}
\end{eqnarray}
The derivatives $\partial_\theta\phi_k(\theta)$ and
$\partial_\theta v_k(\theta)$ exist when $K_\theta(t)$ is differentiable
and can be found by the following algorithm, which is obtained by
differentiating the Durbin--Levinson algorithm:
\begin{eqnarray*}
\partial_\theta\phi_{i,i}(\theta) &=& \Biggl[ \Biggl(
\partial_\theta K_\theta(i\Delta) - \sum_{j=1}^{i-1}
\bigl( \partial_\theta\phi_{(i-1),j}(\theta) K_\theta\bigl((i-j) \Delta\bigr)\\
&&\hspace*{83pt}{} +
\phi_{(i-1),j}(\theta) \partial_\theta K_\theta\bigl((i-j)\Delta\bigr) \bigr) \Biggr)
v_{i-1}(\theta) \\
&&\quad{} + \Biggl( K_\theta(i \Delta) - \sum_{j=1}^{i-1}
\phi_{(i-1),j}(\theta) K_\theta\bigl((i-j) \Delta\bigr) \Biggr) \partial_\theta
v_{i-1}(\theta) \Biggr] v_{i-1}(\theta)^{-2},\\
\pmatrix{\partial_{\theta_j} \phi_{i,1}(\theta)\cr \vdots\cr \partial_{\theta_j} \phi_{i,i-1}(\theta)}
&=&\pmatrix{\partial_{\theta_j}
\phi_{i-1,1}(\theta)\cr
\vdots\cr \partial_{\theta_j} \phi_{i-1,i-1}(\theta)}-
\partial_{\theta_j} \phi_{i,i}(\theta)\pmatrix{\phi_{i-1,i-1}(\theta)
\cr \vdots\cr \phi_{i-1,1}(\theta)}\\
&&{}- \phi_{i,i}(\theta)\pmatrix{\partial_{\theta_j} \phi_{i-1,i-1}(\theta) \cr
\vdots\cr \partial_{\theta_j} \phi_{i-1,1}(\theta)}
\end{eqnarray*}
for $j=1, \ldots, p$, and
\[
\partial_\theta v_i(\theta) = \partial_\theta v_{i-1}(\theta)
\bigl( 1 -\phi_{i,i}(\theta)^2 \bigr) -2 v_{i-1}(\theta)\phi_{i,i}(\theta)\,
\partial_\theta\phi_{i,i}(\theta).
\]

The minimum mean squared error linear predictors of $X((i+1) \Delta)$
and $( X((i+1) \Delta) - \phi_k(\theta)^T X_{i:i+1-k})^2$ given
$X_{i:i+1-k}$ are $\phi_k(\theta)^T X_{i:i+1-k}$ and $v_k(\theta)$,
respectively. This is because for the Gaussian processes considered in
this paper, the two conditional expectations are linear in
$X_{i:i+1-k}$. Thus the pseudo-score function is a prediction-based
estimating function as defined in {S\o rensen}~\cite{ms2011}, where estimating
functions of a slightly more general type than in the original paper
({S\o rensen}~\cite{ms2000}) are treated. The generalization allows the
predicted
function to depend both on the parameter and on the previous
observations. Exploring the
relation of the pseudo-score function to the optimal estimating
function based on these predictors is of interest.

We start by defining a class of prediction-based estimating
functions. Define the $(k+1)\times2$-matrices
\[
Z^{(i)} = \pmatrix{X_{i:i+1-k}^T & 0\vspace*{+2pt}\cr
0 \cdots 0  & 1}^T,\qquad i = k, \ldots, n-1,
\]
and the $(k+1)$-dimensional vectors
\[
H_i(\theta) = Z^{(i)} \pmatrix{X\bigl((i+1) \Delta\bigr) -
\phi_k(\theta)^T X_{i:i+1-k}\cr
\bigl( X\bigl((i+1) \Delta\bigr) - \phi_k(\theta)^T
 X_{i:i+1-k}\bigr)^2 - v_k(\theta)},\qquad i = k, \ldots, n-1.
\]
Then the class of prediction-based estimating functions to which
(\ref{psscore}) belongs is given by
%
%e16 ###
%
\begin{equation}
\label{G}
G_n(\theta) = A(\theta) \sum_{i=k}^{n-1} H_i(\theta),\vspace*{-2pt}
\end{equation}
where $A(\theta)$ is a $p \times(k+1)$ matrix of weights that can depend
on the parameter, but not on the data. The pseudo-score function
(\ref{psscore}) is obtained if the weight matrix
$A(\theta)$ is chosen as
\[
\tilde{A}(\theta) =\pmatrix{\dfrac{\partial_{\theta}
\phi_k(\theta)^T}{v_k(\theta)} & \dfrac{\partial_{\theta}
v_k(\theta)}
{2 v_k(\theta)^2} }
.\vspace*{-1pt}
\]
Within the class of estimators
obtained by solving the estimating equation $G_n(\theta) = 0$ for some
choice of $A(\theta)$, the estimator with the smallest asymptotic variance
is obtained by choosing the optimal weight matrix $A^*(\theta)$. The optimal
estimating function is the one closest to the true score function
in an $L^2$-sense (for details, see Heyde~\cite{heyde97}).

We now can find the optimal weight matrix $A^*(\theta)$. The covariance
matrix of the $(k+1)$-dimensional random vector
$\sum_{i=k}^{n-1} H^{(i)}(\theta)/\sqrt{n-k}$ is
%
%e17 ###
%
\begin{equation}
\label{M}
\bar{M}_n(\theta) = M^{(1)}(\theta) + M^{(2)}_n(\theta),\vspace*{-2pt}
\end{equation}
where
\[
M^{(2)}_n(\theta) = \sum_{j=1}^{n-k-1} \frac{(n-k-j)}{(n-k)}
[ E_\theta(H_k(\theta)H_{k+j}(\theta)^T) +
E_\theta(H_{k+j}(\theta) H_k(\theta)^T ]\vspace*{-2pt}
\]
and
\[
M^{(1)}(\theta) = E_\theta( H_k(\theta) H_k(\theta)^T )
=\pmatrix{v_k(\theta) \mathcal{K}_k(\theta)
& O_{k,1}\cr
O_{1,k} & 2 v_k(\theta)^2},
\]
with $O_{j_1,j_2}$ denoting here and later the $j_1 \times j_2$-matrix of
0s, and with $ \mathcal{K}_k(\theta)$ denoting the covariance matrix
of $(X(k \Delta), \ldots, X(\Delta))$ defined in Section~\ref{sec2}. We have
used $E_\theta(H_i(\theta))=0$, which is a general property of
prediction-based estimating functions (see {S\o rensen}~\cite{ms2011}).
In this
particular case, this is easily seen by finding the conditional
expectation of $H_i(\theta)$ given $X_{i:i+1-k}$, which is~0.

To find the optimal estimating function, we also need the
$p \times(k+1)$ \textit{sensitivity-matrix} $S(\theta)$, given by
%
%e18 ###
%
\begin{equation}
\label{S}
S(\theta)^T = E_\theta(\partial_{\theta^T} H_i(\theta)) =
-\pmatrix{\mathcal{K}_k (\theta) \partial_{\theta^T}
\phi_k(\theta) \cr
\partial_{\theta^T} v_k(\theta)}.\vspace*{-1pt}
\end{equation}
For the derivation of $M^{(1)}(\theta)$ and $S(\theta)$, we use that the model is
Gaussian and, in particular, we use that $\phi_k(\theta)^T
X_{i:i+1-k}$ is the conditional expectation of $X((i+1)\Delta)$ and
not just the minimum
mean squared linear predictor as in the general theory of
prediction-based estimating functions. The optimal
weight matrix is given by
%
%e19 ###
%
\begin{equation}
\label{Astar0}
A^*_n(\theta) = - S(\theta) \bar{M}_n(\theta)^{-1},\vspace*{-2pt}
\end{equation}
see {S\o rensen}~\cite{ms2011}.\vadjust{\goodbreak}

The class of estimating functions considered above is not the
full class of prediction-based estimating functions to which
(\ref{psscore}) belongs, as defined by {S\o rensen}~\cite{ms2000} and
{S\o rensen}~\cite{ms2011}. The full class is obtained by replacing in
(\ref{G})
$A(\theta)$ with a $p \times
2(k+1)$ matrix and $H_i(\theta)$ with the $2(k+1)$-dimensional vectors
$\breve{H}_i(\theta)$ obtained when $Z^{(i)}$ is replaced by the
$2(k+1) \times2$ matrix,
\[
\breve{Z}^{(i)} =\pmatrix{X_{i:i+1-k}^T & 0 & 1 & O_{1,k} \cr
O_{1,k} & 1 & 0 & X_{i:i+1-k}^T}^T
\]
in the definition of $H_i(\theta)$. In this way, $H_i(\theta)$ is
extended by $k+1$ extra coordinates. Because the moments of an odd
order of a centered multivariate Gaussian distribution equal 0, we see that
the extra $k+1$ coordinates of $\breve{H}_i(\theta)$ have expectation
0 under the true probability measure irrespective of the value of
the parameter $\theta$; therefore,
they cannot be expected to be a useful addition to $H_i(\theta)$.
However, the
extra coordinates might be correlated with the
coordinates of $H_i(\theta)$, and thus might be used to reduce the
variance of the estimating function. To see that this is not the case,
the optimal estimating function based on $\breve{H}_i(\theta)$ can be
calculated.\vspace*{1pt} The covariance matrix of the random vector
$\sum_{i=k}^{n-1} \breve{H}^{(i)}(\theta)/\sqrt{n-k}$ can be shown to
be a block-diagonal matrix with two $(k+1) \times(k+1)$-blocks, the
first of which equals $\bar{M}_n(\theta)$. This follows from the
fact that moments of an odd order of a centered multivariate Gaussian
distribution equal 0. Moreover, the sensitivity matrix
corresponding to $\breve{H}_i(\theta)$ is
\[
\breve{S}(\theta)^T = E_\theta(\partial_{\theta^T} \breve
{H}_i(\theta)) =
- \pmatrix{\mathcal{K}_k (\theta)\,
\partial_{\theta^T} \phi_k(\theta) \cr
\partial_{\theta^T} v_k(\theta)\cr
O_{k+1,p}}.
\]
Therefore, the optimal weight matrix is
\[
\breve{A}^*_n(\theta) = \pmatrix{ A^*_n(\theta) &
O_{p,k+1}       },
\]
and thus the optimal prediction-based estimating function obtained
from $\breve{H}_i(\theta)$ equals the optimal estimating function obtained
from $H_i(\theta)$. It is therefore sufficient to consider the
aforementioned smaller class of prediction-based estimating functions,
which we do in the rest of the paper.

The pseudo-score function, $\partial_\theta\tilde{\ell}_n(\theta)$,
is not equal to the optimal prediction-based estimating function.
In fact,
%
%e20 ###
%
\begin{equation}
\label{Atilde}
\tilde{A}(\theta) = - S(\theta) M^{(1)}(\theta)^{-1}.
\end{equation}
The magnitude of the difference between the two estimating functions
depends on how small the entries of $M^{(2)}_n(\theta)$
are relative to the entries of $M^{(1)}(\theta)$. Because correlations
decrease exponentially with the distance in time (see {Rei\ss}~\cite{reiss},
page 26), the terms in the sum
defining $M^{(2)}_n(\theta)$ can be small compared with the entries of
$M^{(1)}(\theta)$; however, under what conditions this occurs and exactly
how small the terms are depend on $\theta$, $\Delta$, and $k$.

In the next section we show that the limit
%
%e21 ###
%
\begin{equation}
\label{M2}
M^{(2)}(\theta) = \lim_{n \rightarrow\infty} M^{(2)}_n(\theta) =
\sum_{j=1}^\infty[ E_\theta(H_k(\theta)H_{k+j}(\theta)^T) +
E_\theta(H_{k+j}(\theta) H_k(\theta)^T )]
\end{equation}
exists. Therefore, we can define the following weight matrix, which
does not
depend on~$n$:
%
%e22 ###
%
\begin{equation}
\label{Astar}
A^*(\theta) = - S(\theta) \bar{M}(\theta)^{-1},
\end{equation}
where
%
%e23 ###
%
\begin{equation}
\label{Mb}
\bar{M}(\theta) = \lim_{n \rightarrow\infty} \bar{M}_n(\theta) =
M^{(1)}(\theta) + M^{(2)}(\theta).
\end{equation}
The estimating function
%
%e24 ###
%
\begin{equation}
\label{Gstar}
G^*_n(\theta) = A^*(\theta) \sum_{i=k}^{n-1} H_i(\theta)
\end{equation}
is asymptotically optimal and theoretically is easier to handle than
$A_n^*(\theta) \sum_{i=k}^{n-1} H_i(\theta)$. In practice, the
optimal weight matrices $A_n^*(\theta)$ or $A^*(\theta)$ usually must
be calculated by simulation. The amount of computation can be reduced
by using the
approximation to $G^*_n(\theta)$ obtained by replacing $A^*(\theta)$ or
$A_n^*(\theta)$ with the matrix obtained from (\ref{Astar}) and
(\ref{Mb}) when $M^{(2)}(\theta)$ is replaced by a suitably truncated
version of the series in (\ref{M2}). This does not make much
difference, because the terms in the sum (\ref{M2}) decrease
exponentially fast.
%
%s4 ###
\section{Asymptotics of the pseudo-likelihood estimator}\label{sec4}
In this section, we present the asymptotic properties of estimators obtained
by solving the estimating equation $G_n(\hat{\theta}_n) = 0$, where
$G_n$ is given by (\ref{G}). Important particular cases of this are
the maximum
pseudo-likelihood estimator obtained by maximizing (\ref{pslik}) and the
optimal prediction-based estimator obtained by solving $G_n^*(\hat
{\theta}_n)
= 0$ with $G_n^*$ given by (\ref{Gstar}). The depth, $k$, of $G_n$ is
assumed fixed.
The asymptotic properties are proven for a solution to the general equation
(\ref{geneqn}) under the following assumption:

\begin{con}
\label{conAA}
\begin{longlist}[(b)]
\item[(a)] The functions $K_\theta(t)$ and $A(\theta)$ are twice
continuously differentiable with respect to~$\theta$.
\item[(b)] The $p \times(k+1)$ matrix
$(\partial_\theta\phi^T_k(\theta)\enskip
\partial_\theta v_k(\theta))$ has rank $p$ (in particular,
$k+1 \geq p$).
\item[(c)] $A(\theta)\bar\mathcal{K} (\bar{\phi}_k(\theta_0) -
\bar{\phi}_k(\theta))=0$ if and only if $\theta= \theta_0$.
\end{longlist}
\end{con}

Here
%
%e25 ###
%
\begin{equation}
\label{Kbar}
\bar\mathcal{K} = \pmatrix{\mathcal{K}_k(\theta_0) & O_{k,1}\cr
O_{1,k} & 1},
\end{equation}
and
\[
\bar{\phi}_k(\theta) = \pmatrix{\phi_k(\theta)\cr
v_k(\theta) +
2 \phi_k(\theta)^T \kappa_k (\theta_0) - \phi_k(\theta)^T
\mathcal{K}_k(\theta_0)\phi_k(\theta) }.
\]

If $A$ equals $\tilde{A}$ (corresponding to the pseudo-score
function) or $A$ equals $A^*$ (corresponding to the
optimal prediction-based estimating function), then Condition
\ref{conAA}(a) is satisfied if $K_\theta(t)$ is three times continuously
differentiable, which is the case for the models considered in Examples
\ref{affinex2} and~\ref{affinex4}. Condition~\ref{conAA}(a) ensures
that the functions $\phi_k(\theta)$ and $v_k(\theta)$ are continuously
differentiable. Condition~\ref{conAA}(c) is an identifiability
condition that ensures eventual uniqueness of the estimator.
\begin{theo}
\label{thA1}
Assume that the true parameter value $\theta_0$ belongs to the
interior of the parameter space $\Theta$. Suppose that Condition
\ref{conAA} is satisfied, and that the matrix $A(\theta_0)$ has full
rank $p$. Then a consistent estimator $\hat{\theta}_n$ that solves the
estimating equation $G_n(\hat{\theta}_n)=0$ exists and is unique in
any compact subset of $\Theta$ containing $\theta_0$ with a
probability tending to 1 as $n \rightarrow\infty$. Moreover,
\[
\sqrt{n}(\hat{\theta}_n- \theta_0) \stackrel{\mathcal
{D}}{\longrightarrow}
N_p ( 0,U(\theta_0)^{-1}V(\theta_0) (U(\theta_0)^{-1})^T )
\]
as $n \rightarrow\infty$, where $V(\theta_0) = A(\theta_0) \bar{M}
(\theta_0) A(\theta_0)^T$ with $\bar{M}(\theta_0)$ given by
(\ref{Mb}), and
%
%e26 ###
%
\begin{equation}
\label{U}
U(\theta_0)
=E_{\theta_0} (\partial_{\theta} G_n(\theta_0)^T)/(n-k) = S(\theta_0)
A(\theta_0)^T.
\end{equation}
Here $S(\theta)$ is the sensitivity matrix given by (\ref{S}).
\end{theo}

Note that it follows from (\ref{Astar0}) and (\ref{Atilde}) that
$A^*(\theta_0)$ and $\tilde{A}(\theta_0)$ have rank $p$ if
$S(\theta_0)$ has rank $p$, because $\bar{M}_n$ and $M^{(1)}$ are
non-singular covariance matrices. That $S(\theta_0)$ has rank $p$
follows from Condition~\ref{conAA}(b) by (\ref{S0}) below.
\begin{pf*}{Proof of Theorem~\ref{thA1}} The theorem follows from general
asymptotic statistical
results for stochastic processes (see, e.g., Jacod and {S\o rensen}
\cite{jjms}). We need to
establish that a law of large numbers and a central limit theorem
hold and to check regularity conditions.

Under our general assumption that $X$ is stationary, {Rei\ss} \cite
{reiss} (page
25) showed that $X$ is exponentially $\beta$-mixing. Therefore,
a law of large numbers holds for sums of the form $n^{-1}\sum_{i=1}^n
f(X_{i+k:i})$. The process $\{ H_i(\theta_0) \}$ is exponentially
$\alpha$-mixing, and because the process $X$ is Gaussian, $H_i(\theta)$
has moments of all orders. Therefore, it follows from Theorem 1 in
Section~1.5 of Doukhan~\cite{doukhan} that (\ref{M2}) converges, and that
\[
\frac{G_n(\theta_0)}{\sqrt{n}} \stackrel{\mathcal
{D}}{\longrightarrow}
N ( 0,V(\theta_0) )
\]
as $n \rightarrow\infty$.

Next, we need to check regularity conditions that ensure the
asymptotic results. The estimating function satisfies
that $E_{\theta_0}(G_n(\theta_0)) = 0$. Furthermore, it is obvious
from the
continuous differentiability of the functions $\phi_k(\theta)$ and
$v_k(\theta)$ that the derivatives
$\partial_{\theta_j} A(\theta) H_i(\theta) + A(\theta) \,\partial
_{\theta_j}
H_i(\theta)$ are locally dominated\vadjust{\goodbreak} integrable. Finally, the matrix
$U(\theta_0)$ is invertible because $A(\theta_0)$ has full rank and
%
%e27 ###
%
\begin{equation}
\label{S0}
S(\theta_0) = - \pmatrix{\partial_\theta\phi^T_k(\theta_0)& \partial
_\theta v_k(\theta_0) } \bar\mathcal{K},
\end{equation}
with $\bar\mathcal{K}$ given by (\ref{Kbar}). The first matrix has full
rank by Condition~\ref{conAA}(b), and $\bar\mathcal{K}$
is invertible because $\mathcal{K}_k(\theta_0)$ is the covariance matrix
of $X(\Delta), \ldots, X(k\Delta)$, which is not degenerate. Now the
existence and consistency of $\hat{\theta}_n$, as well as the eventual
uniqueness of a consistent estimator on any compact subset of $\Theta$
containing $\theta_0$, follow (see Jacod and {S\o rensen}~\cite{jjms}).
The locally dominated
integrability of $A(\theta)H_i(\theta)$ (which follows from Condition
\ref{conAA}(a)) implies that
$n^{-1}G_n(\theta)$ converges uniformly to
$A(\theta)E_{\theta_0}(H_i(\theta)) = A(\theta)\bar\mathcal{K}
(\bar{\phi}_k(\theta_0) -
\bar{\phi}_k(\theta))$ for $\theta$ in a compact set.
The fact that the limit is a continuous functions of $\theta$ and
satisfies $A(\theta)E_{\theta_0}(H_i(\theta)) \neq0$ for $\theta
\neq\theta_0$ implies that any non-consistent solution to the
estimating equation will eventually leave any compact subset of
$\Theta$ containing $\theta_0$. The asymptotic normality follows by
standard arguments (see, e.g., Jacod and {S\o rensen}~\cite{jjms}).
\end{pf*}

A simpler estimator with the same asymptotic distribution as in the
estimator from (\ref{G}) is obtained from the estimating function
\[
G^\circ_n(\theta) = A(\theta^\circ_n) \sum_{i=k}^{n-1} H_i(\theta),
\]
where $\theta^\circ_n$ is some consistent estimator of $\theta$,
obtained, for
instance, by simply using $p$ suitably chosen coordinates of
$H_i(\theta)$. For this estimating function, the identifiability
condition Condition~\ref{conAA}(c) can be replaced by the following
condition:
\begin{con}
\label{conAAA}
\begin{longlist}[(b)]
\item[(a)] The function $(\phi_k(\theta),v_k(\theta))$ is one-to-one.
\item[(b)] $\bar{\phi}_k(\theta_0) - \bar{\phi}_k(\theta) \in N^\bot$
for all
$\theta\in\Theta$, where $N$ is the null space of the matrix
$A(\theta_0)\bar\mathcal{K}$.
\end{longlist}
\end{con}

This readily follows from the fact that the limit of $n^{-1}G^\circ
_n(\theta)$
is $ A(\theta_0)\bar\mathcal{K} (\bar{\phi}_k(\theta_0) -
\bar{\phi}_k(\theta))$. In the case of the pseudo-likelihood function,
we have the simple expression
\[
\tilde{A}(\theta_0)\bar\mathcal{K} = \pmatrix{ v_k(\theta_0)^{-1}\,
\partial_\theta\phi_k^T(\theta_0) \mathcal{K}_k(\theta_0)&
\tfrac12 v_k(\theta_0)^{-2}\,
 \partial_\theta v_k(\theta_0)
}.
\]
Condition~\ref{conAAA}(a) is a basic assumption without which there
is no hope of estimating $\theta$ using the pseudo-score function
(\ref{psscore}). The condition must be checked for individual
models. Obviously, it is not always satisfied, as demonstrated by the
following examples. Consider the model in Example~\ref{affinex2} with
the restriction that $b=a$. For this model, the autocovariance function
depends on $r$ and $b$ only through $r-b^{-1}>0$ for $t \in[0,r]$. In
Example~\ref{affinex4} with the restriction that $a=0$, the
autocovariance function depends on $r$ and $b$ only through $r \sqrt{b}$
for $t \in[0,r]$.

Theorem~\ref{thA1} implies that the asymptotic distribution
of the optimal prediction-based estimator, $\hat{\theta}^*_n$, is
%
%e28 ###
%
\begin{equation}
\label{asyopt}
\sqrt{n}(\hat{\theta}^*_n- \theta_0) \stackrel{\mathcal
{D}}{\longrightarrow}
N_p ( 0, ( S(\theta_0) \bar{M}(\theta_0)^{-1} S(\theta_0)^T
)^{-1} ),\vadjust{\goodbreak}
\end{equation}
and the asymptotic distribution of the pseudo-likelihood estimator,
$\tilde{\theta}_n$, is
%
%e29 ###
%
\begin{equation}
\label{asypseu}
\sqrt{n}(\tilde{\theta}_n- \theta_0) \stackrel{\mathcal
{D}}{\longrightarrow}
N_p \bigl( 0, W(\theta_0)^{-1} + W(\theta_0)^{-1}B(\theta_0)
W(\theta_0)^{-1}\bigr ),
\end{equation}
where
\[
W(\theta) = S(\theta) M^{(1)}(\theta)^{-1} S(\theta)^T =
\frac{ \partial_\theta{\phi}_k(\theta)^T \mathcal{K}_k(\theta)\,
\partial_{\theta^T} {\phi}_k(\theta)}{v_k(\theta)} +
\frac{\partial_\theta v_k(\theta) \,\partial_{\theta^T} v_k(\theta)}
{2 v_k(\theta)^2}
\]
and
\[
B(\theta) = \tilde{A}(\theta)M^{(2)}(\theta)\tilde{A}(\theta)^T
= S(\theta) M^{(1)}(\theta)^{-1} M^{(2)}(\theta)
M^{(1)}(\theta)^{-1} S(\theta)^T.
\]
The result for $\hat{\theta}^*_n$ follows because
\[
- S(\theta_0) A^*(\theta_0)^T = A^*(\theta_0) \bar{M}
(\theta_0) A^*(\theta_0)^T = S(\theta_0) \bar{M}(\theta)^{-1}
S(\theta_0)^T,
\]
and the result for $\tilde{\theta}_n$ follows because
\[
- S(\theta_0) \tilde{A}(\theta_0)^T = S(\theta_0) M^{(1)}(\theta_0)^{-1}
S(\theta_0)^T
\]
and
\[
\tilde{A}(\theta_0) \bar{M} (\theta_0) \tilde{A}(\theta_0)^T =
S(\theta_0) M^{(1)}(\theta_0)^{-1} S(\theta_0)^T +
\tilde{A}(\theta_0) M^{(2)}(\theta_0) \tilde{A}(\theta_0)^T.
\]

According to the general theory of estimating functions (see, e.g.,
Heyde~\cite{heyde97}), the matrix $S(\theta_0) \bar{M}(\theta_0)^{-1}
S(\theta_0)^T
-(W(\theta_0)^{-1} + W(\theta_0)^{-1}B(\theta_0)W(\theta_0)^{-1})^{-1}$
is positive definite; that is, the asymptotic covariance matrix of
$\tilde{\theta}_n$ is larger than that of $\hat{\theta}^*_n$ (in the
usual ordering of positive semi-definite matrices). Thus
the asymptotic variance of $f(\tilde{\theta}_n)$ is larger than
that of $f(\hat{\theta}^*_n)$ for any differentiable function $f\dvtx \R^p
\mapsto\R$. If $B(\theta_0)$ is invertible, then
\[
[W(\theta_0)^{-1} + W(\theta_0)^{-1}B(\theta_0)W(\theta_0)^{-1}]^{-1}
=W(\theta_0) - [B(\theta_0)^{-1} + W(\theta_0)^{-1}]^{-1},
\]
and if $M^{(2)}(\theta_0)$ is invertible, then
\[
\bar{M}(\theta_0)^{-1} = M^{(1)}(\theta_0)^{-1} -
M^{(1)}(\theta_0)^{-1} \bigl[M^{(1)}(\theta_0)^{-1} +
M^{(2)}(\theta_0)^{-1}\bigr]^{-1}M^{(1)}(\theta_0)^{-1},
\]
where we have used twice that $(I+A)^{-1} = I - A(I+A)^{-1}$ for a
matrix $A$.
Thus the difference between the two inverse asymptotic covariance
matrices can be expressed as
%
%e30 ###
%
\begin{eqnarray}\label{loss}
&&S(\theta_0) \bar{M}(\theta_0)^{-1} S(\theta_0)^T -
[W(\theta_0)^{-1} + W(\theta_0)^{-1}B(\theta_0)W(\theta
_0)^{-1}]^{-1}\nonumber\\
&&\quad=[B(\theta_0)^{-1} + W(\theta_0)^{-1}]^{-1} \nonumber\\
&&\qquad{}-
S(\theta_0)M^{(1)}(\theta_0)^{-1} \bigl[M^{(1)}(\theta_0)^{-1} +
M^{(2)}(\theta_0)^{-1}\bigr]^{-1}M^{(1)}(\theta_0)^{-1}S(\theta_0)^T \\
&&\quad=\bigl[ \bigl(\tilde{A}(\theta_0) M^{(1)}(\theta_0)
\tilde{A}(\theta_0)^T \bigr)^{-1} + \bigl(\tilde{A}(\theta_0)M^{(2)}(\theta_0)
\tilde{A}(\theta_0)^T \bigr)^{-1} \bigr] ^{-1} \nonumber\\
&&\qquad{}- \tilde{A}(\theta_0) \bigl[M^{(1)}(\theta_0)^{-1} +
M^{(2)}(\theta_0)^{-1} \bigr]^{-1} \tilde{A}(\theta_0)^T.\nonumber
\end{eqnarray}
%
%This is an expression of the superiority of the optimal prediction-based
%estimator over the pseudo-likelihood estimator.

It is considerably easier to calculate the pseudo-likelihood function
(\ref{pslik}) than the optimal estimating function (\ref{Gstar}),
because the latter involves derivatives with respect to $\theta$ of
the covariance function and higher-order moments of $X$. In
particular, in cases where the covariance function is not explicitly
known and must be determined by simulation, it is much easier to
calculate (\ref{pslik}) than~(\ref{Gstar}). Thus the maximum pseudo-likelihood
estimator is preferred in practice. The formula~(\ref{loss}) then can
be used to assess whether
the loss of efficiency relative to the optimal estimator is
acceptable.
\begin{ex}
\label{lossex}
As an example, we calculated the efficiency loss for the model (\ref{affin})
in Example~\ref{affinex2} in a number of cases. When $k$
is sufficiently large, the pseudo-likelihood function is almost
efficient, and thus the information loss (\ref{loss}) is necessarily
small. Therefore, it is most interesting to calculate the efficiency
loss when $k$ is small. We calculate the relative information loss,
that is, the information loss (\ref{loss}) relative to the information
for the optimal estimator given by (\ref{asyopt}). The main problem is to
calculate the matrix (\ref{M2}). However, for $k=1$, a simple
expression for each term in the sum (\ref{M2}) can be obtained using
the formula of Isserlis~\cite{isserlis}, and so a suitably truncated
version of
the sum can be easily calculated. For $k \geq2$, the matrix (\ref{M2})
can be determined by simulation using (\ref{Mb}).
Specifically, we determine the covariance matrix
$\bar M_n(\theta_0)$, with $n$ suitably large, by simulation. This is
computationally more demanding.

We first considered the efficiency loss for the parameter $b$ in the
case $k=1$. The parameters $\sigma^2$ and $r$ were fixed at a value
of 1, and $a$ was chosen equal to $-1$. For $b = - \mathrm{e}^{-2} = -0.1353$ (the
value for which the mixing rate is maximal), the relative information
loss was found to be very small, less than 0.1 percent for $\Delta= 0.1$, $\Delta
= 0.5$, and $\Delta=1$.

Next, we calculated the information loss for a number of values of $b$
with $\Delta=1$. For $b = -0.3, -0.06, 0.05, 0.1, 0.2, 0.3$, and $0.5$,
the mixing rate is relatively high, and the information loss is
less than 0.1 percent. For $b = -0.5, -0.4, 0.7$, and $0.9$, the information loss is
between 0.1 and 1 percent, whereas for $b=-0.6, -0.7$, and $-0.9$, it is
1.4 percent, 3.0 percent, and 9.8 percent, respectively.

Finally, we calculated the relative the information loss for $k=3$ and
$k=5$. In this case, information loss was calculated for both $a$
and $b$. The parameters $a$, $\sigma^2$, and $r$ had the same value as
before, and $\Delta=1$. For $b = - \mathrm{e}^{-2}$, the relative information
loss for both $a$ and $b$ is  less than 0.1 percent for both values
of $k$. For $b = -0.5$ the information loss is less than 0.1 percent for $b$
and 0.2 percent for $a$.

In most cases, the relative information loss is so tiny that in
practice it is preferable to use the maximum pseudo-likelihood
estimator. The information loss increases as the mixing rate
decreases. Only for $k=1$ and $b = -0.9$ is the information loss large
enough to justify the use of the more complicated optimal estimator.
\end{ex}
%
%s5 ###
\section{Simulation study}
\label{simstudy}
In this section we report the results of a simulation study in which we
investigated some
properties of the pseudo-likelihood estimator introduced in Section
\ref{sec2}. We restrict ourselves to the model considered in Example
\ref{affinex2} and to estimating $\theta= (a,b)$.
The delay time $r$ is chosen equal to 1, and $\sigma^2$ is
fixed at 1. This study was not intended to serve as a complete
simulation study; rather, the intention was to illustrate some
properties of the estimator
and give a first impression of how the joint distribution
of the two-dimensional estimator $\tilde{\theta}_n = (\tilde{a}_n,
\tilde{b}_n)$ depends on the time between observations $\Delta$,
the depth $k$ of the pseudo-likelihood function, and the true
parameter value. We performed simulations for three values of
$\theta$: $\theta= (-1,0.95)$ near the upper boundary of the domain of
stationarity, $\theta= (-1, -1/\mathrm{e}^2) = (-1, -0.1353)$ which is the
parameter value with the highest possible mixing rate for the
stationary solution $X$ when $a=-1$, and $\theta= (-1,-2.1)$
near the lower boundary of the domain of stationarity.
For each parameter value, we considered four sampling frequencies
with the same number of observation time points (200); specifically, the
observation time points were $i \Delta$, $i=1, \ldots, 200$, with
$\Delta= 0.05, 0.1, 0.5, 1$. The simulations of the SDDE were
conducted with a step size of 0.001. In all cases, 1000 data sets were
simulated, and thus 1000 estimates were generated. For each data set, a
new trajectory of the driving Wiener process was generated. The full
simulation study is reported in {K\"uchler} and {S\o rensen}~\cite{kuso1}.

Table~\ref{table4b} reports
the mean values and standard deviations of the simulated estimates of
$a$ and $b$ for $(a,b) = (-1,-0.1353)$. For $(a,b) = (-1,0.95)$, the
estimators are more biased and have a larger standard deviation for
small values of $\Delta$ and $k$, whereas for $(a,b) = (-1,-2.1)$, the
bias is small and the standard deviations are comparable in all cases. For
$(a,b) = (-1,0.95)$, the estimators of $a$ and $b$ are highly
correlated, whereas this is the case only for small values of $\Delta$
and $k$ for the other parameter values.%

%
%t1 ###
%
\begin{table}
\caption{Mean and standard deviation of the pseudo-likelihood
estimator of $a$ (upper part of the table)
and $b$ (lower part of the table) for various values of
depth $k$, time
between observations~$\Delta$, and number of observations $n$.
In all
cases $n \Delta$, the length of the observation interval, is 200, and
the true parameter values are $a=-1$ and $b=-0.1353$}
\label{table4b}
\begin{tabular*}{\textwidth}{@{\extracolsep{\fill}}d{1.2}d{4.0}d{2.2}
d{2.2}d{2.2}d{2.2}d{2.2}d{2.2}d{2.2}@{}}
  \hline
   &  & \multicolumn{7}{l@{}}{\hspace*{-2pt}$k$} \\     [-5pt]
   &  & \multicolumn{7}{c@{}}{\hspace*{-2pt}\hrulefill}\\
  \multicolumn{1}{@{}l}{$\Delta$} & \multicolumn{1}{l}{$n$} &
  \multicolumn{1}{l}{\hspace*{-2pt}1} &
  \multicolumn{1}{l}{\hspace*{-2pt}3} & \multicolumn{1}{l}{\hspace*{-2pt}5}
  & \multicolumn{1}{l}{\hspace*{-2pt}7} &
  \multicolumn{1}{l}{\hspace*{-2pt}9} & \multicolumn{1}{l}{\hspace*{-2pt}13} &
  \multicolumn{1}{l@{}}{\hspace*{-2pt}20} \\ \hline
0.05 & 4000 & -1.73 & -1.07 & -1.04 & -1.02 & -1.02 & -1.01 & -1.01 \\
& & 2.32 & 0.21 & 0.14 & 0.11 & 0.11 & 0.10 & 0.09 \\[3pt]
0.1 & 2000 & -1.27 & -1.03 & -1.03 & -1.01 & -1.01 & -1.01 & -1.01 \\
& & 0.83 & 0.13 & 0.10 & 0.09 & 0.09 & 0.09 & 0.09 \\ [3pt]
0.5 & 400 & -1.04 & -1.01 & -1.01 & -1.01 & -1.01 & -1.01 & -1.01 \\
& & 0.14 & 0.10 & 0.09 & 0.09 & 0.10 & 0.09 & 0.10 \\ [3pt]
1.0 & 200 & -1.02 & -1.01 & -1.01 & -1.02 & -1.02 & -1.01 & -1.02 \\
& & 0.12 & 0.11 & 0.11 & 0.11 & 0.11 & 0.11 & 0.12 \\[3pt]
2.0 & 100 & -1.10 & -1.04 & -1.03 & -1.03 & -1.04 & -1.04 & -1.04 \\
& & 0.43 & 0.21 & 0.19 & 0.18 & 0.21 & 0.19 & 0.17 \\[15pt]
0.05 & 4000 & 0.42 & -0.09 & -0.11 & -0.15 & -0.13 & -0.14 & -0.14 \\
& & 2.66 & 0.44 & 0.31 & 0.23 & 0.19 & 0.14 & 0.09 \\[3pt]
0.1 & 2000 & 0.08 & -0.13 & -0.14 & -0.14 & -0.14 & -0.13 & -0.13 \\
& & 1.14 & 0.27 & 0.18 & 0.13 & 0.11 & 0.09 & 0.09 \\[3pt]
0.5 & 400 & -0.12 & -0.14 & -0.14 & -0.13 & -0.14 & -0.14 & -0.14 \\
& & 0.28 & 0.10 & 0.11 & 0.11 & 0.11 & 0.11 & 0.10 \\[3pt]
1.0 & 200 & -0.14 & -0.14 & -0.14 & -0.13 & -0.14 & -0.13 & -0.13 \\
& & 0.16 & 0.13 & 0.13 & 0.13 & 0.13 & 0.13 & 0.14 \\[3pt]
2.0 & 100 & -0.26 & -0.15 & -0.14 & -0.13 & -0.15 & -0.14 & -0.14 \\
& & 0.65 & 0.29 & 0.27 & 0.25 & 0.30 & 0.27 & 0.24  \\  \hline
\end{tabular*}
\end{table}

The most remarkable observations from our simulation study can be
summarized as follows:
\begin{itemize}
\item For a fixed number of observation time points, the
bias and standard deviation of the estimators worsen as the time
between observations $\Delta$ decreases, at least when $\Delta\leq
r$. For $\Delta> r$, the quality does not change much with $\Delta$,
and whether the bias and variance
increase or decrease with $\Delta$ depends on the parameter value.
\item The smaller the $\Delta$ value, the more the choice of the
depth $k$ of the pseudo-likelihood functions influences the quality of
the estimators when $\Delta\leq r$. For $\Delta> r$, the importance
of $k$ increases again for some parameter values.
\item It is surprising that a similar pattern is seen when the length
of the observation interval $n \Delta$ is fixed so that the sample
size decreases as $\Delta$ increases. However, here there is a
clearer tendency for the estimators to deteriorate when $\Delta> r$,
so that there is an optimal value of~$\Delta$, which seems to be
around $r$.
\item The absolute value of the correlation between
$\tilde{a}$ and $\tilde{b}$ decreases with increasing depth $k$ to a
limit, which is strongly dependent on the true parameter value. Near the
upper boundary of the stability region, the estimators are highly
correlated. A high absolute value of the correlation indicates that
it is difficult to distinguish between the effects of the lagged term and
the nonlagged term in the drift; thus, it is not surprising that the
absolute correlation is large when the depth is small.
\item For small values of the depth $k$, the joint distribution of the
estimators of $a$ and $b$ can deviate from a two-dimensional normal
distribution by having crescent-shaped contours.
\end{itemize}
\section*{Acknowledgements}
We thank Adrian M. K\"ammler for assisting with the efficiency loss
calculations, Katja Krol for writing the PC program used in the
simulation study, and Daniel Skodlerack for providing an earlier
version of that program. We also thank the referees, the associate
editor, and the
editor for their constructive criticism that led to considerable
improvement of the manuscript. The research of Michael S\o rensen was
supported by the Danish Center for Accounting and Finance, funded by
the Danish Social Science Research Council; by the Center for
Research in Econometric Analysis of Time Series, funded by the Danish
National Research Foundation; and by a grant from the University of
Copenhagen Excellence Programme.
%
% imsref loaded by audrone.aklyte, 2012-02-22 13:09:33
% imsref loaded by audrone.aklyte, 2012-02-22 13:25:25
%

\printhistory

\end{document}